\documentclass[11pt]{elsart}
\usepackage{amssymb,amsmath}
\usepackage{mathrsfs}
\usepackage{epic}
\newtheorem{theo}{Theorem}[section]
\newtheorem{lema}{Lemma}[section]

\newcommand{\racion}[2]{\mbox{\small$\frac{{#1}}{{#2}}$}}

\newcommand\be{\begin{enumerate}}
\newcommand\ee{\end{enumerate}}
\newcommand\bi{\begin{itemize}}
\newcommand\ei{\end{itemize}}

\newcommand{\pe}[2]{\left\langle#1,#2\right\rangle}

\newcommand\XX[1]{\mathbb{#1}}

\newcommand\bd{\begin{definicion}{\bf }}
\newcommand\ed{\end{definicion}}
\newcommand\bl{\begin{lema}{\bf }}
\newcommand\el{\end{lema}}
\newcommand\bp{\begin{prop}{\bf }}
\newcommand\ep{\end{prop}}
\newcommand\bt{\begin{theo}{\bf }}
\newcommand\et{\end{theo}}
\newcommand\bdm{\begin{proof}}
\newcommand\edm{\end{proof}}
\newcommand\bn{\begin{nota}{\bf }}
\newcommand\en{\end{nota}}
\newcommand\bc{\begin{corolary}{\bf }}
\newcommand\ec{\end{corolary}}

\newtheorem{nota}{Remark}[section]

\newtheorem{definicion}{Definition}[section]
\newtheorem{corolary}{Corolary}[section]
\newenvironment{proof}{{\bf Proof:\ }}{\hfill
$\square$}

\newfont{\got}{eufm10 scaled \magstep1}

\usepackage{bm,hyperref}

\begin{document}
\begin{frontmatter}
\title{The complementary polynomials
and the Rodrigues operator.
A distributional study${}^\star$}
\thanks{This work has been supported by
Direcci´on General de Investigaci´on
(Ministerio de Educaci´on y Ciencia) of
Spain, grant MTM 2006-13000-C03-02.}
\date{\today}
\author{Roberto S. Costas-Santos}
\address{Department of Mathematics,
University of California, South Hall,
Room 6607, Santa Barbara, CA 93106, USA}
\ead{rscosa@gmail.com}
\ead[url]{http://www.rscosa.com}
\begin{abstract}
We can write the polynomial solution of
the second order linear differential
equation of hypergeometric-type
$$
\phi(x)y''+\psi(x)y'+\lambda y=0,
$$
where $\phi$ and $\psi$ are polynomials,
$\deg \phi\le 2$, $\deg \psi=1$ and
$\lambda$ is a constant, among others,
by using the Rodrigues operator
$R_k(\phi,{\bf u})$ (see \cite{coma2})
where $\bf u$ is certain linear operator
which satisfies the distributional equation
\begin{equation} \label{1}
\frac{d}{dx}[\phi {\bf u}]=\psi {\bf u},
\end{equation}
as
$$
P_n(x)=B_n R_n(\phi,{\bf u})[1],
\qquad B_n\ne 0,\quad n=0, 1, 2, \dots
$$
Taking this into account we construct the
complementary polynomials.
Among the key results is a generating functional
function in closed form leading to derivations
of recursion relations and addition theorem.
The complementary polynomials satisfy a
hypergeometric-type differential equation
themselves, have a three-term recursion
among others and Rodrigues formulas.
\end{abstract}
\end{frontmatter}

\section{Introduction}

In this paper\footnote{Note that this paper
is a natural extension of \cite{web1}
in terms of linear functionals.} we
will consider the po\-ly\-no\-mi\-al solutions,
$P_n$, of the second order linear differential
equation of hypergeometric-type \cite{niuv1}
\begin{equation} \label{2}
\phi(x)P''_n(x)+\psi(x)P'_n(x)+\lambda_n P_n(x)=0,
\end{equation}
where $\phi$ and $\psi$ are polynomials, $\deg \phi
\le 2$, $\deg \psi\le 1$, and
$$
\
\lambda_n=-n\psi'-\frac{n}{2}(n-1)\phi'',
\quad n= 0,1, 2, \dots
$$
In fact, it is well-known that such family of
polynomials satisfies a property of orthogonality,
i.e.
$$
\pe {{\bf u}}{P_n P_m}=r_n \delta_{n,m},\quad
r_n\ne 0, \quad n,m= 0, 1, 2, \dots
$$
where $\delta_{n,m}$ is the Kronecker's symbol
and the linear functional $\bf u$ satisfies
the distributional equation
\begin{equation} \label{3}
\frac{d}{dx}[\phi {\bf u}]=\psi {\bf u}.
\end{equation}
Taking this into account is possible
to writhe such polynomials trough
the following Rodrigues functional formula:
\begin{equation} \label{4}
P_n {\bf u}=B_n \frac{d^n}{dx^n}[{\bf u}_n], \quad
B_n\ne 0, \quad n= 0, 1, 2, \dots
\end{equation}
where ${\bf u}_0:={\bf u}$, ${\bf u}_k:=\phi
{\bf u}_{k-1}$, $k= 1, 2, 3, \dots$.
\bn
Note that although is possible to find an
integral representation for the linear
functional $\bf u$, since the polynomial
solutions of such differential equation are
the classical orthogonal polynomials (COP)
\cite{niuv1}, is important to take into
account that the technique presented in this
manuscript can be generalized to the
semiclassical orthogonal polynomials where
such integral representation is not known
for each semiclassical functional (see, for
example, \cite{med2,aratma}).
\en

Our first goal is to construct complementary
polynomials for them by using their Rodrigues
functional representation, Eq. \eqref{4}, in a
simple and natural way.
The generating functional function of these
complementary polynomials is obtained in
closed form.

The paper is organized as follows.
In the next section we introduce
some notations and definitions useful
for the next sections. In section 3
we construct the complementary polynomials.
In Section 4 we establish their generating
functional function.
The second order linear differential
equation associated to the
complementary polynomials is derived in
Section 5.

\section{Preliminaries}
In this section we will give a brief
survey of the operational calculus that
we will use in the rest of the paper.

Let $\XX P$ be the linear space of
polynomial functions in $\XX C$ (in the
following we will refer to them as
polynomials) with complex coefficients and
$\XX P'$ be its algebraic dual space, i.e.,
$\XX P'$ is the linear space of all linear
applications ${\bf u}:\XX P \to  \XX C$.
In the following we will call the elements
of $\XX P'$ as functionals.

Let $\{P_n\}_{n\ge 0}$ be a sequence of
polynomials such that $\deg P\le n$ for all
$n\ge 0$.
A sequence defined in this way is said to
be a basis or a basis sequence of $\XX P$
if and only if $deg P_n=n$ for all $n\ge 0$.
Since the elements of $\XX P'$ are linear
functionals, it is possible to determine
them from their actions on a given basis
$\{P_n\}_{n\ge 0}$ of $\XX P$.

In general, we will represent the action of a
functional over a polynomial by
$$
\pe {{\bf u}} P, \quad {\bf u}\in \XX P',
\quad P\in \XX P.
$$
Therefore, a functional is completely determined
by a sequence of complex numbers $u_n:=\pe
{{\bf u}} {x^n}$, $n\ge 0$, the so-called moments
of the functional.

\bd
Let ${\bf u}\in \XX P'$ be a functional.
We say that $\bf u$ is a quasi-definite
(or regular) functional if and only if
there exists a polynomial sequence $\{P_n
\}_{n \ge 0}$, which is orthogonal with
respect to $\bf u$.
\ed


Let us define the following operations in
$\mathbb{P}'$.
For any polynomial $h$ and any
$c\in \mathbb{C}$, let ${\bf u}':=\frac
d{dx}{\bf u}$, $h{\bf u}$, and $(x-c)^{-1}
{\bf u}$ be the linear functionals defined
on $\mathbb{P}$ by
(see \cite{mar3,khe})
\be
\item $\pe {{\bf u}'} {P}:=-\pe {{\bf u}}{P'},
\quad P\in\mathbb{P}$,
\item $\pe {Q{\bf u}} P :=\pe {{\bf u}}{Q\, P},
\quad P,\;Q\in\mathbb{P},$
\item $\pe {(x-c)^{-1}{\bf u}} P:=\pe {{\bf u}}
{\theta_c(P)}, \quad P\in\mathbb{P},\quad
\mbox{where}\;\theta_c(P)(x)=\frac{P(x)-P(c)}{x-c}.$
\ee
Furthermore, for any linear functional
$\bf u$ and any polynomial $P$
we get
\begin{equation} \label{relgu}
\frac{d}{dx}[P{\bf u}]:=(P{\bf u})'=
P{\bf u}'+g'{\bf u}.
\end{equation}

\section{Complementary Polynomials}

Before to introduce such polynomials let
us introduce the Rodrigues operator.
\bd \label{rod-op-def}
Given a polynomial $\phi$ and a linear
functional ${\bf u}$, such that there
exists a polynomial $\psi$ in such a way
$$
\frac{d}{dx}[\phi {\bf u}]=\psi {\bf u}.
$$
We define the {\bf $k$-th Rodrigues
operator} associated with the pair
($\phi$, $\bf u$), as:
$$
\begin{array}{rl}
R_0(\phi,{\bf u}):=& I, \\
R_1(\phi,{\bf u})[p]=& q \quad
where \quad \dfrac{d}{dx}[p\,
{\bf u}_1]=q\, {\bf u},\\ R_k(
\phi,{\bf u}):=& R_1(\phi,{\bf
u})\circ R_{k-1}(\phi,{\bf u}_1),
\quad k=2,3,\dots
\end{array}
$$
where $I$ represents the identity
operator.
\\
The pair $(\phi, {\bf u})$ satisfying
the above condition we will call
classical pair.
\ed
The following results are related
with the Rodrigues operator.
\bl \label{lem-latt}
\cite[Lemma 4.2]{coma2}
Let $(\phi, \bf u)$ be a classical
pair, then for any positive integer
$k$ and any polynomial $\pi$, the
function
$$
R_1(\phi,{\bf u}_k)[\pi]
$$
is a polynomial of degree $\deg(\pi)+1$
with leading coefficient:
$$
-\frac{\lambda_{m+2k}}{m+2k}\ne 0,
\quad k= 0, 1, 2, \dots, \quad m=
\deg \pi.
$$
\el
\bn
In the sequel we will denote by
$\mathscr{R}_{k,\ell}$ the $k$-th
Rodrigues operator associated with
the pair ($\phi$, ${\bf u}_k$),
$R_k(\phi,{\bf u}_\ell)$, and
$\mathscr{R}_{k}:=\mathscr{R}_{k,1}$.
\en
Let us now introduce the complementary
polynomials, $\mathscr{P}_\nu(x;n)$,
defining them in terms of the
Rodrigues operator.
\bd
Let $n, \nu$ be two integers, with
$0\le \nu \le n$, the complementary
polynomial of degree $\nu$ with
respect to the $P_n$ is defined as
\begin{equation} \label{5}
P_n(x)=B_n\mathscr{R}_{n-\nu}[
\mathscr{P}_\nu(x;n)],\quad B_n\ne 0.
\end{equation}
\ed
\bn
Note these polynomials for $\nu=1$
are considered by M. Alfaro and R.
\'{A}lvarez-Nodarse in \cite{alal}.
More generally, is easy to check that
$\mathscr{P}_{n-\nu}(x;n)$ is, up to
a constant, the $\nu$-th derivative of
$P_n$, i.e. $P^{(\nu)}_n(x)$ (see e.g.
\cite[Lemma 3.3]{coma2}).
\en
Taking into account this remark and the
theory of COP, the following result
holds.
\bl
${\mathscr P}_\nu(x;n)$ is a polynomial
of degree $\nu$ that satisfies the
recursive differential equation:
\begin{equation} \label{6}
\mathscr{P}_{\nu+1}(x;n)=\phi(x)\frac{d
\mathscr{P}_\nu(x;n)}{dx}+\big(\psi(x)+
(n-\nu-1)\phi'(x)\big)\mathscr{P}_\nu(x;n).
\end{equation}
By the Rodrigues functional formula
\eqref{4}, $\mathscr{P}_{0}(x;n)\equiv 1.$
\el
\bdm
Equations \eqref{5}; and  \eqref{6}
follow by induction.
Case $\nu=1$ is derived by carrying out
explicitly the innermost differentiation
in Eq. \eqref{4}, which is a natural way
of working with the Rodrigues functional
formula \eqref{4} that yields
\begin{equation}%
P_n(x){\bf u}=B_n\frac{d^{n-1}}{dx^{n-1}}
\frac{d}{dx}[\phi^n {\bf u}]=B_n\frac{d^{n
-1}}{dx^{n-1}}\left[n\phi^{n-1}\phi'{\bf
u}+\phi^n {\bf u}'\right].
\label{6.1}
\end{equation}%
showing, taking into account \eqref{3}),
that
$$
P_n(x)=B_n \mathscr{R}_{n-1}[(n-1)
\phi'(x)+\psi(x)],
$$
therefore
\begin{equation} \label{7}
\mathscr{P}_1(x;n)=(n-1)\phi'(x)+\psi(x).
\end{equation}
Assuming the validity of the Rodrigues
functional formula \eqref{5} for $\nu$
we carry out another differentiation in
Eq. \eqref{5} obtaining
\begin{eqnarray}\nonumber
P_n(x){\bf u}&=&B_n\frac{d^{n-\nu-1}}{dx^{n-
\nu-1}}\big\{(n-\nu)\phi'\mathscr{P}_\nu
(x;n)\phi^{n-\nu-1}{\bf u}\\ \nonumber&+&
\phi^{n-\nu}(x)\mathscr{P}_\nu(x;n){\bf u}'+
\phi^{n-\nu}\mathscr{P}'_\nu(x;n){\bf u}
\big\}\\\nonumber&=&B_n\frac{d^{n-
\nu-1}}{dx^{n-\nu-1}}\big[\big\{(n-\nu)\phi'
\mathscr{P}_{\nu}(x;n)
(\psi-\phi')\mathscr{P}_\nu(x;n)+\phi
\mathscr{P}'_\nu(x;n)\big\}\phi^{n-\nu-1}
{\bf u}\big]\\ \label{7.1} &=&B_n\frac{d^{n-
\nu-1}}{dx^{n-\nu-1}}\left(\mathscr{P}_{\nu
+1}(x;n){\bf u}_{n-\nu-1}\right).
\end{eqnarray}
Comparing the right hand side of Eq. \eqref{7}
proves Eq. \eqref{5} by induction along with
the recursive differential equation (DE)
\eqref{6}.
\edm

In terms of a generalized Rodrigues
functional representation we have
\bt \label{theo2.1}
The polynomials $\mathscr{P}_{\nu}
(x;n)$ satisfy the Rodrigues functional
formulas
\begin{equation} \label{8}
\mathscr{P}_{\nu}(x;n){\bf u}_{n-\nu}=
\frac{d^{\nu}}{dx^{\nu}}[{\bf u}_n]
\quad \iff \quad \mathscr{P}_{\nu}(x;n)=
\mathscr{R}_{\nu,n-\nu}[1],
\end{equation}
\begin{equation} \label{9}
\mathscr{P}_{\nu}(x;n){\bf u}_{n-\nu}=\frac
{d^{\nu-\mu}}{dx^{\nu-\mu}}[\mathscr{P}_{\mu}
(x;n){\bf u}_{n-\mu}]\quad \iff \quad
\mathscr{P}_{\nu}(x;n)=\mathscr{R}_{\nu-\mu,
n-\nu}\circ \mathscr{R}_{\mu,n-\mu}[1].
\end{equation}
\et
Note that the proof of this result is
straightforward taking into account the
hypergeometric character of the DE \eqref{2}
and Lemma 3.3 in \cite{coma2}.

\section{generating functional Function}
\bd  The generating functional function for the
polynomials $\mathscr{P}_\nu(x;n)$ is
\begin{equation} \label{9.1}
\mathscr{P}(y,x;n)=\sum_{\nu=0}^{\infty}
\frac{y^{\nu}}{\nu!}\mathscr{P}_{\nu}(x;n).
\end{equation}
The series converges for $|y|<\epsilon$
for some $\epsilon>0$ and can be summed
in closed form if the generating functional
function is regular at the point $x$.
\ed
\bt
The generating functional function
for the polynomials $\mathscr{P}_\nu
(x;n)$ is given in closed form by
\begin{equation} \label{10}
\mathscr{P}(y,x;n){\bf u}=\left(\frac{\phi(x+
y\phi(x))}{\phi(x)}\right)^n \widetilde{\bf u},
\end{equation}
\begin{equation} \label{11.1}
\frac{\partial^\mu}{\partial y^\mu}
\mathscr{P}(y,x;n){\bf u}=\left(
\frac{\phi(x+y\phi(x))}{\phi(x)}
\right)^{n-\mu} \mathscr{P}_\mu(x+y
\phi(x);n)\widetilde{\bf u},
\end{equation}
where $\widetilde{\bf u}=\exp(y\phi(x)
\racion d{dx})\circ {\bf u}$.
\et
\bdm Taking into account that ${\bf u}_k=
\phi^k {\bf u}$, equation \eqref{10}
follows by multiplying the generatriz
function by $\phi^n$ and substituting
the Rodrigues functional representation,
Eq. \eqref{8} in Eq. \eqref{9.1} which
yields
$$
\phi^n(x)\mathscr{P}(y,x;n){\bf u}=
\mathscr{P}(y,x;n){\bf u}_n=
\sum_{\nu=0}^{\infty}\frac{y^{\nu}}{\nu!}
\mathscr{P}_{\nu}(x;n)\phi^n(x){\bf u}=
\sum_{\nu=0}^{\infty}\frac{(y\phi(x))^{\nu}}
{\nu!}\frac{d^{\nu}}{dx^{\nu}}[{\bf u}_n],
$$
converging for $|y\phi(x)|<\epsilon$ for a
suitable $\epsilon>0$ if $x$ is a regular
point of the generating functional function.

In fact, at this point the series can be
summed exactly, because the expression
inside the derivatives is independent
of the summation, obtaining
\begin{equation}
\label{11}
\phi^n(x)\mathscr{P}(y,x;n){\bf u}=
\exp(y\phi(x)\racion{d}{dx})[\phi^n(x){\bf u}]=
\phi^n(x+y\phi(x))\widetilde {\bf u},
\end{equation}
and therefore
\begin{equation} \label{12}
\mathscr{P}(y,x;n){\bf u}=\left(\frac{
\phi(x+y\phi(x))}{\phi(x)}\right)^n
\widetilde {\bf u}.
\end{equation}

Differentiating Eq. \eqref{9.1} and
substituting  the generalized Rodrigues
functional formula \eqref{9} in this
yields Eq. \eqref{11.1} similarly.
\edm
\bn
Note that if the linear $\bf u$
admits an integral representation
as
$$
\pe {{\bf u}}{P}=\int_\Omega P(z)\rho(z)
dz,
$$
where $\Omega$ is certain contour in the
complex plane, then we can rewrite \eqref{12}
as
\begin{equation} \label{12.2}
\mathscr{P}(y,x;n)=\left(\frac{
\phi(x+y\phi(x))}{\phi(x)}\right)^n
\frac{\rho(x+y\phi(x)}{\rho(x)}.
\end{equation}
\en
Since we are going to consider recurrence
 relations we are going to translate the
case $\mu=1$ of Eq. \eqref{11.1} into partial
differential equations (PDEs).
\\
\bt
The generating functional function
satisfies the PDEs
\begin{equation} \label{12.1}
(1+y\phi'(x)+\racion{1}{2}y^2\phi''
\phi(x))\frac{\partial \mathscr{P}
(y,x;n)}{\partial y}=[\mathscr{P}_1
(x;n)+y\phi(x)\mathscr{P}_1'(x;n)]
\mathscr{P}(y,x;n),
\end{equation}
\begin{equation} \label{13}
\frac{\partial \mathscr{P}(y,x;n)}
{\partial y}=[(n-1)\phi'(x+y\phi(x))
+\psi(x+y\phi(x))]\mathscr{P}(y,x;n-1),
\end{equation}
\begin{eqnarray}\nonumber
&&\left(1+y\phi'(x)+\racion{1}{2}y^2\phi''
\phi(x)\right)\frac{\partial \mathscr{P}
(y,x;n)}{\partial x}\\ \label{14}&=&\bigg\{
(1+y\phi'(x))\mathscr{P}_1'(x;n)-\frac{1}{2}
y\phi''\mathscr{P}_1(x;n)
\bigg\}y\mathscr{P}(y,x;n),
\end{eqnarray}
\begin{eqnarray}\nonumber
\phi(x)\frac{\partial \mathscr{P}
(y,x;n)}{\partial x}&=&(1+y\phi'(x)
)[\psi(x)+(n-1)\phi'(x)+y\phi(x)
(\psi'+(n-1)\phi'')]\\&\times&
\mathscr{P}(y,x;n-1)-[\psi(x)+(n-
1)\phi'(x)]\mathscr{P}(y,x;n).
\label{15}
\end{eqnarray}
\et
\bdm
From Eq. \eqref{11.1} for $\mu=1$ in
conjunction with Eq. \eqref{10} and
taking into account that $\bf u$ is
regular we obtain
\begin{equation} \label{16}
\phi(x+y\phi(x))\frac{\partial
\mathscr{P}(y,x;n)}
{\partial y}=\phi(x)[\psi(x+y\phi(x))
+(n-1)\phi'(x+y\phi(x))]\mathscr{P}
(y,x;n).
\end{equation}
Substituting in Eq. \eqref{16} the
Taylor series-type expansions
\begin{eqnarray}\nonumber
\phi(x+y\phi(x))&=&\phi(x)(1+y\phi'(x)+
\racion{1}{2}y^2\phi''\phi(x)),\\\nonumber
\phi'(x+y\phi(x))&=&\phi'(x)+y\phi''\phi(x),\\
\psi(x+y\phi(x))&=&\psi(x)+y\psi'(x)\phi(x),
\label{17}
\end{eqnarray}
since $\phi$ and $\psi$ are polynomials of
degree, at most, 2 and 1, respectively, we
verify Eq. \eqref{12.1}.
Using the exponent $n-1$ instead
of $n$ of the generating functional
function we can obtain Eq. \eqref{13}.

By differentiation of the
generating functional function, Eq.
\eqref{11}, with respect
to the variable $x$ we find
Eq. \eqref{15}.
Using the exponent $n$ instead
of $n-1$ of the generating
functional function in conjunction
with Eq. \eqref{15} leads to Eq.
\eqref{14} where in any case
we were assume the linear
functional is regular.
\edm

In fact the polynomials
$\mathscr{P}_\nu(x;n)$
satisfy different recursions
and other relations
which are the same that
H. J. Weber obtained in \cite{web1}
hence we will omite here.

\section{The second order linear differential
equation of hypergeometric-type}

\bt
The polynomials $\mathscr{P}_{\nu}(x;n)$
satisfy the Sturm-Liouville differential equation
\begin{equation} \label{18}
\frac{d}{dx}\left(\frac{d\mathscr{P}_\nu(x;n)}
{dx}{\bf u}_{n-\nu+1}\right)=
-\mu_{n,\nu} \mathscr{P}_\nu(x;n)
{\bf u}_{n-\nu},
\end{equation}
which is equivalent to, assuming $\bf u$ is
regular,
\begin{equation} \label{20}
\phi(x)\frac{d^2\mathscr{P}_\nu(x;n)}
{dx^2}+[(n-\nu)\phi'(x)+\psi(x)]
\frac{d\mathscr{P}_\nu(x;n)}{dx}=
-\mu_{n,\nu} \mathscr{P}_\nu(x;n),
\end{equation}
and the eigenvalues are given by
\begin{equation} \label{19}
\mu_{n,\nu}=-\nu ((n-\racion{\nu+1}{2})
\phi''+\psi'), \qquad \nu=0, 1, \ldots.
\end{equation}
\et
Taking into account the hypergeometric
character of the DE \eqref{1}, and Theorem
\ref{theo2.1} this result is very well-known
and straightforward.

\section{Conclusions}

We have used a natural way of working with
the Rodrigues operator and the Rodrigues
functional formula of a given set polynomials
solutions of the second order differential
equation \eqref{1} which leads to a set of
closely related complementary polynomials
that satisfy their own Rodrigues formulas,
always have a generating functional function
that can be summed in closed form leading to
numerous recursion relations and addition
theorems.

Furthermore, since the differential equation
\eqref{1} is of hypergeometric-type, the
complementary polynomials satisfy a
second order linear differential equation
similar to that of the original polynomials.

In fact, it is a straightforward calculation
to verify that this method generates all the
basics of the COP, $\Delta$-classical OP, and
$q$-polynomials (see, for example,
\cite{coma2},\cite{mealma},\cite{niuv3} and
references therein).

Moreover this method, considered in a similar
way by H. J. Weber, is not restricted to
the COP and can be applied to semiclassical
orthogonal polynomials.

\bibliographystyle{plain}

\begin{thebibliography}{99}
\bibitem{alal}
M. Alfaro and R. \'{A}lvarez-Nodarse.
\newblock A characterization of the classical orthogonal discrete and
  $q$-polynomials.
\newblock {\em J. Comput. Appl. Math.}  {\bf 201} (2007), 48--54.

\bibitem{aratma}
J. Arves\'{u}, J. Atia, and F. Marcell\'{a}n.
\newblock Semiclassical linear functionals: The symmetric companion.
\newblock {\em Communications in the Analytic Theory of Continued Fractions},
  {\bf 10} (2002), 13--29.

\bibitem{coma2}
R. S. Costas-Santos and F. Marcellan.
\newblock Classical orthogonal polynomials. A general difference calculus
  approach.
\newblock 2008. Submitted

\bibitem{khe}
L. Kheriji.
\newblock An introduction to the $H_q$-semiclassical orthogonal polynomials.
\newblock {\em Meth. Appl. Anal.} {\bf  10} (2003), 387--412.

\bibitem{mar3}
P. Maroni.
\newblock Semiclassical character and finite-type relations between polynomial
  sequences.
\newblock {\em Appl. Num. Math.} {\bf  31} (1999), 295--330.

\bibitem{med2}
J. C. Medem.
\newblock A family of singular semi-classical functionals.
\newblock {\em Indag. Mathem.} {\bf  13}(3) (2002), 351--362.

\bibitem{mealma}
J. C. Medem, R. \'{A}lvarez-Nodarse, and F. Marcell\'{a}n.
\newblock On the $q$-polynomials: a distributional study.
\newblock {\em J. Comput. Appl. Math.} {\bf  135} (2001), 157--196.

\bibitem{niuv3}
A. F. Nikiforov and V. B. Uvarov.
\newblock Classical orthogonal polynomials in a discrete variable on nonuniform
lattices.
\newblock {\em Inst. Prikl. Mat. Im. M. V. Keldysha Akad. Nauk SSSR}, 17, 1983.

\bibitem{niuv1}
A. F. Nikiforov and V. B. Uvarov.
\newblock {\em The Special Functions of Mathematical Physics}.
\newblock Birkh\"{a}user Verlag, Basel, 1988.

\bibitem{web1}
H. J. Weber.
\newblock Connections between real polynomial solutions of hypergeometric-type
  differential equations with rodrigues formula.
\newblock {\em Cent. Eur. J. Math.} {\bf  5}(2) (2007), 415--427.
\end{thebibliography}

\end{document}